\documentclass{elsarticle}

\usepackage{lineno,hyperref}
\modulolinenumbers[5]

\journal{Journal of Computational Physics}

\usepackage{amssymb}
\usepackage{amssymb}
\usepackage{amsthm}
\usepackage{lineno}

\usepackage[inline]{asymptote}

\usepackage{float}
\usepackage{multirow}
\usepackage{booktabs}
\usepackage[dvipsnames]{xcolor} % for setting colors
\usepackage{array}
\usepackage{xspace}   	   % Pour avoir des espaces malins dans les formules
\usepackage{mathrsfs}
\usepackage{mdwlist}         % Pour des listes sans trop d'espaces.
\usepackage{placeins}
\usepackage[english]{babel}
\usepackage{mathtools}
\usepackage[linesnumbered,lined,boxed,commentsnumbered]{algorithm2e}

\usepackage{bigints}
\usepackage{tensor}
\usepackage{multicol}
\usepackage{subfig}
\usepackage{tikz}
\usepackage{tikz-3dplot}
\usetikzlibrary{calc}
\usetikzlibrary{patterns}
\usetikzlibrary{arrows}
\usetikzlibrary{arrows.meta}
\usetikzlibrary{shadings}
\usetikzlibrary{positioning}
\usepackage{pgfplots}
\usepackage{listings}

\usepackage[normalem]{ulem}

\usepackage{dsfont}

\newcommand{\inputtikzfile}[1]{%
\IfFileExists{./tikz/#1.pdf}{\includegraphics[scale=1]{./tikz/#1.pdf}}{\input{tikz/#1.tikz}}%
}

\makeatletter
\DeclareFontFamily{U}{tipa}{}
\DeclareFontShape{U}{tipa}{m}{n}{<->tipa10}{}
\newcommand{\arc@char}{{\usefont{U}{tipa}{m}{n}\symbol{62}}}%

\newcommand{\arc}[1]{\mathpalette\arc@arc{#1}}

\newcommand{\arc@arc}[2]{%
  \sbox0{$\m@th#1#2$}%
  \vbox{
    \hbox{\resizebox{\wd0}{\height}{\arc@char}}
    \nointerlineskip
    \box0
  }%
}
\makeatother

\definecolor{liquid}{rgb}{0.28,0.46,1}

\makeatletter
\newcommand\resetsubfigs{\setcounter{sub\@captype}{0}}
\makeatother

\begin{document}

\begin{frontmatter}

\title{Quaternionic octahedral fields\\ $\text{SU}(2)$ parameterization of 3D frames}

\author[a,b]{Pierre-Alexandre Beaufort\corref{mycorrespondingauthor}}
\author[a]{Jonathan Lambrechts}
\author[b]{Christophe Geuzaine}
\author[a]{Jean-Fran\c cois Remacle}

\address[a]{Universit\'e catholique de Louvain, iMMC, Avenue Georges Lemaitre 4, 1348 Louvain-la-Neuve, Belgium}
\address[b]{Universit\'e de Li\`ege, Montefiore Institute, All\'ee de la D\'ecouverte 10, B-4000 Li\`ege, Belgium}%\fntext[myfootnote]{Since 1880.}

\cortext[mycorrespondingauthor]{Corresponding author}

\begin{abstract}
   Full hexahedral meshes are required to be as regular as possible, which means that the local topology has to be constant almost everywhere.
  This constraint is usually modelled by 3D frames.
  A 3D frame consists of three mutual orthogonal (unit) vectors, defining a local basis.
  3D frame fields are auxiliary for hexahedral mesh generation.
  Computation of 3D frame fields is an active research field.
  There mainly exist three ways to represent 3D frames: combination of rotations, spherical harmonics and fourth order tensor.
  We propose here a representation carried out by the special unitary group.
  The article strongly relies on \cite{du1964homographies}.
  We first describe the rotations with quaternions, \cite[\S 13-15]{du1964homographies}.
  We define and show the isomorphism between unit quaternions and the special unitary group, \cite[\S 16]{du1964homographies}.
  The frame field space is identified as the quotient group of rotations by the octahedral group, \cite[\S 20]{du1964homographies}.
  The invariant forms of the vierer, tetrahedral and octahedral groups are successively built, without using homographies \cite[\S 39]{du1964homographies}.
  Modifying the definition of the isomorphism between unit quaternions and the special unitary group allows to use the invariant forms of the octahedral group as a unique parameterization of the orientation of 3D frames.
  The parameterization consists in three complex values, corresponding to a coordinate of a variety which is embedded in a three complex valued dimensional space.
  The underlined variety is the model surface of the octahedral group, which can be expressed with an implicit equation.
  We prove that from a coordinate of the surface, we may identify all the quaternions giving the corresponding 3D frames.
  We show that the euclidean distance between two coordinates does not correspond to the actual distance of the corresponding 3D frames.
  We derive the expression of three components of a coordinate in the case of frames sharing an even direction.
  We then derive a way to ensure that a coordinate corresponds to the special unitary group.
  Finally, the attempted numerical schemes to compute frame fields are given.
  \end{abstract}


\end{frontmatter}

\linenumbers

%%%%%%%%%%%%%%%%%%
\section{Introduction}

Full hexahedral mesh is still an open question (see \citep{shepherd2008hexahedral}).
Yet, it seems that there is an easy way to produce a full hexahedral mesh: first produce a tetrahedral mesh, then split each of them into 4 hexahedra.
But this way is not convenient: the hexahedra are not regular, they tend to have bad quality and do not form a structured mesh.
Finite element community aims to get full hexahedral meshes, possibly structured, which are as regular as possible.	

The	regularity of an hexahedral mesh is related to the topology of a given domain $R\subset \Re^3$.
Let us consider a mesh on $R \subset \Re^3$ with $N$ nodes (i.e. vertices), $N_E$ edges, $N_F$ facets and $N_C$ cells (i.e. element-wise volumes, here being hexahedra) is such that
\begin{equation}\label{eq:euler3d}
\chi(R) = N - N_E + N_F - N_C
\end{equation}
with $\chi(R)$ the Euler characteristic of the region $R$.
The Euler characteristic of a region is half the one of its boundary \citep[§4C, (4-15)]{Gross2004Jun}
\begin{equation}\label{eq:euler2d3d}
\chi(R) = \dfrac{\chi(\partial R)}{2}
\end{equation}
We assume there are $n$ nodes, $n_e$ edges and $n_f$ facets making the mesh of $\partial R$.
\begin{equation}\label{eq:euler3d2d}
\chi(R) = \dfrac{1}{2} \left( n - n_e + n_f \right)
\end{equation}

From a topological point of view, an hexahedral mesh is assumed to be regular if each inner (boundary) vertex is shared by
\begin{itemize}
\item  8 (4) hexahedra
\item  6 (5) edges
\item  12 (8) facets
\end{itemize}
whose corresponding equations are
\begin{equation}\label{eq:regHex}
\left.\begin{array}{rcl}
8 N_V &=& 8(N-n) + 4n\\
2 N_E &=& 6(N-n) + 5n\\
4 N_F &=& 12(N-n) + 8n\\
\end{array}\right\}
\end{equation}
Using \eqref{eq:regHex} into \eqref{eq:euler3d}, we get
$$
8\chi(R) = 0
$$
It means that if a region $R$ may be meshed by regular hexahedra, its Euler characteristic is zero.
But the opposite is not true: a region whose characteristic is zero does not mean that it may be meshed  by regular hexahedra.
For example, let us consider a region that is meshed such that there are $k$ inner loops, each made of $L$ edges.
Those $kL$ edges are shared by  3 hexahedra; those edges are then singular.
The $N-kL$ remaining edges are regular.
Equations \eqref{eq:regHex} become
\begin{equation}\label{eq:innerLoop}
\left.\begin{array}{rcl}
8 N_V &=& 8(N-n-kL) + 4n + 6kL\\
2 N_E &=& 6(N-n-kL) + 5n + 5kL\\
4 N_F &=& 12(N-n-kL) + 8n + 9kL\\
\end{array}\right\}
\end{equation}

\begin{figure}
\begin{center}
\includegraphics[scale=.5]{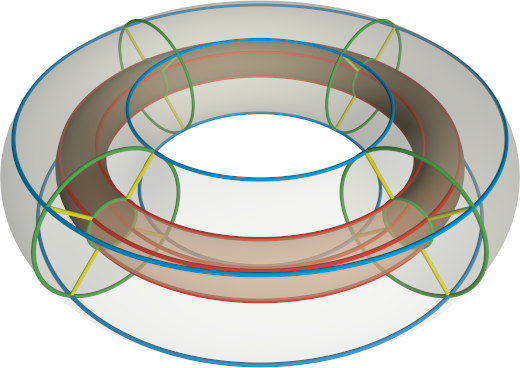}
\caption{Block structure decomposition of a torus and cut torus. The red lines are singular if the inner torus is not an hole, but 4 block structures (front, back, left \& right parts). The front part of the tori is more transparent than the other parts.}
\label{fig:torus}
\end{center}
\end{figure}

Again, using \eqref{eq:innerLoop} into \eqref{eq:euler3d} gives an Euler characteristic that is zero.
While the Euler characteristic defines completely the topology of an oriented 2-manifold (surface), it is not the case for an oriented 3-manifold (region).
Indeed, from \eqref{eq:euler2d3d} a full torus and a torus cut by a smaller one (i.e. the larger one contains the smaller one, Fig. \ref{fig:torus}) have the same Euler characteristic, which is zero. 
Obviously, the cut torus may be meshed by regular hexahedra: you produce a regular quadrangulation of the outer boundary that is mapped onto the inner one, then you link the corresponding vertices.
A full torus cannot be meshed with regular hexahedra; its block structure decomposition corresponds to four singular inner loops.
Both situations are represented by Fig. \ref{fig:torus}.
Unfortunately, topological constraints for hexahedrizations are not as nicely summarized as the ones for quadrangular meshes, \citep[§1, (7)]{Beaufort2017Jan}.

\begin{figure}[!ht]
\begin{center}
\subfloat[Inner vertex connectivity.]{\includegraphics[width=.32\linewidth]{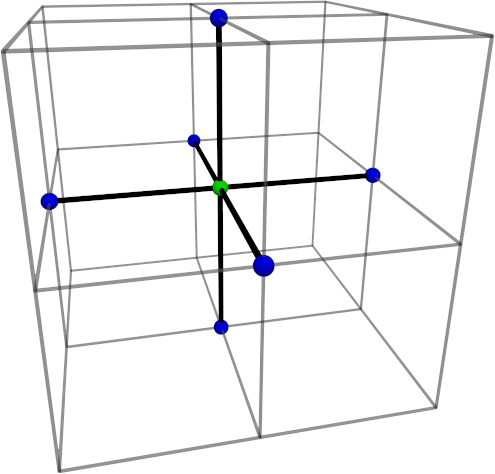}\label{sub::hexgrid}}\hspace*{.25cm}
\subfloat[3D frame.]{\includegraphics[width=.32\linewidth]{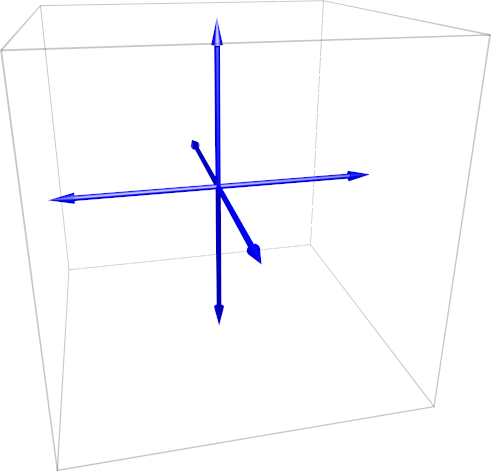}\label{sub::frame}}\hspace*{.25cm}
\subfloat[Octahedron.]{\includegraphics[width=.32\linewidth]{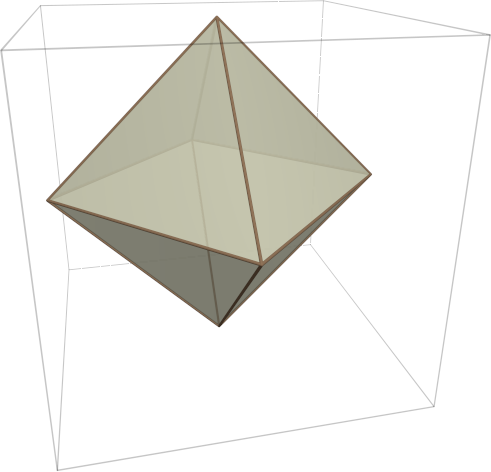}\label{sub::octahedron}}
\caption{Hexahedral features.}
\end{center}
\end{figure}

In order to build full hexahedral mesh that are as regular as possible, we use a three-dimensional frame field designing the desired connectivity of a regular hexahedral mesh, Fig. \ref{sub::hexgrid}.
A 3D frame field gives in each point a 3D frame, picturing the local orientation (and thus the vertex connectivity) of an hexahedron.
Since an orientation is relative, it is measured from the cartesian frame, which is the reference 3D frame Fig. \ref{sub::frame}.
Observe that the corresponding vector field is symmetric, since a frame shares the symmetries of an octahedron, Fig. \ref{sub::octahedron}.

%%%%%%%%%%%%%%%%%%
\section{Rotation Representations}

There exist various ways to represent a frame field.
But at the end of the day, they all essentially consist in rotations of a vector field representing an object that exhibits the 24 symmetries of the octahedral group.
Such objects may be for example fourth order tensors \citep{chemin2018representing}, or spherical harmonics \citep{ray2016practical}.
In those latter cases, they are represented by a nine-dimensional vector.
Actually, they are both based on the representation of
\begin{equation}\label{eq:sphere4}
\hat{f}(x;y;z) = x^4 + y^4 + z^4
\end{equation}
which is the polynomial exhibiting the 24 octahedral symmetries corresponding to the cartesian frame.

% SPHERICAL HARMONICS

\begin{figure}
\begin{center}
\includegraphics[width=.48\linewidth]{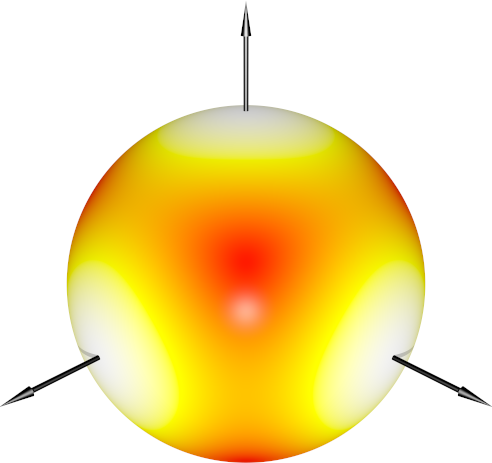}
\caption{$\hat{f}(x;y;z)=x^4+y^4+z^4$ on the sphere $S_2$.}
\label{fig:polynomialOnSphere}
\end{center}
\end{figure}

In the case of \emph{spherical harmonics}, it is seen as a polynomial taking values on the sphere $S_2$, Fig. \ref{fig:polynomialOnSphere}.
This polynomial may be decomposed with the real spherical harmonics of fourth degree.
$$
\hat{f}(x;y;z) = \dfrac{4\sqrt{\pi}}{15}\left( Y_{4,0}(x;y;z) +  \sqrt{\dfrac{5}{7}}~Y_{4,4}(x;y;z)\right)+\dfrac{3}{5}, \forall (x;y;z)\in S_2
$$
with
$$
\begin{array}{rcl}
\left.Y_{4,0}\right|_{S_2} &=& \dfrac{3}{16}\sqrt{\dfrac{1}{\pi}} \left( 3x^4 + 3y^4 + 8z^4 + 6x^2y^2 - 24x^2z^2 - 24y^2z^2 \right)\\
\left.Y_{4,4}\right|_{S_2} &=& \dfrac{3}{16}\sqrt{\dfrac{35}{\pi}} \left( x^4 + y^4 - 6x^2y^2 + 12\right)
\end{array}
$$
If we rotate the cartesian frame $\hat{f}$ by means of a matrix $R\in \text{SO}(3)$, we get $f$\footnote{by denoting $x=x_0,~y=x_1,~z = x_2$}
$$
f(x_0;x_1;x_2) := \hat{f}(R_{0i}x_i;R_{1j}x_j;R_{2k}x_k)
$$
with $x_m = R_{mn} x_n$.
The function $f$ may still be expressed with real spherical harmonics of fourth degree.
% FOURTH ORDER TENSORS

If we consider the isosurface described by the points where polynomial \eqref{eq:sphere4} is equal to one, it corresponds to a unit sphere in 4-norm which may be written as a \emph{fourth order tensor}
\begin{figure}
\begin{center}
\includegraphics[scale=.5]{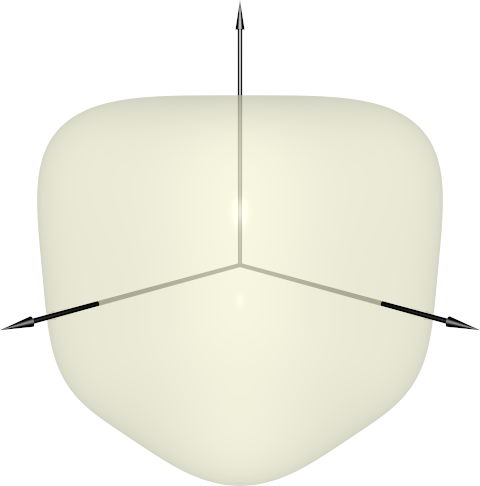}
\caption{$\hat{f}(x;y;z)=1$.}
\end{center}
\end{figure}
$$
\begin{array}{rcl}
\hat{f}(x_0;x_1;x_2) &=& \hat{A}_{ijkl} x_i x_j x_k x_l \\
\text{with}~ \hat{A}_{ijkl} &=& \displaystyle{\sum_{q=1}^3 \delta_{iq} \delta_{jq} \delta_{kq} \delta_{lq}}
\end{array}
$$
Again, if we rotate the frame $\hat{f}$ with a matrix $R\in \text{SO}(3)$, we have
$$
f(x_0,x_1,x_2) = \hat{f}\left(R_{1i}x_i,R_{2i}x_i,R_{3i}x_i\right)
$$
which generalizes $\hat{A}$ as fourth order tensor
$$
A_{ijkl} = R_{im} R_{jn} R_{ko} R_{lp} \hat{A}_{mnop}
$$
giving
$$
f(x_0,x_1,x_2) = A_{ijkl} x_i x_j x_k x_l
$$

As illustrated by those two representations, we understand that 3D frame fields consist in rotations depicted by the quotient group
$$
\text{SO}(3) / \text{O} \subset \text{SO}(3)
$$
where $\text{O}$ is the octahedral group, i.e. the 24 rotations leaving invariant the orientation of an octahedron whose vertices are at the units of each axis, Fig. \ref{sub::octahedron}.

Those two representations work with an object sharing the octahedral symmetries, which enables to identify this quotient group.
We here propose to work directly with the corresponding rotational group.
To do so, we need to describe three rotational groups by means of quaternions.
We later build the corresponding invariant forms by avoiding their symmetries.

\subsection{Quaternions}

A real quaternion $q$ consists of four real numbers $(q_0;q_1;q_2;q_3)\in\Re^4$.
Using three imaginary units $\mathbf{i},\mathbf{j},\mathbf{k}$ such that
\begin{itemize}
\item $\mathbf{i}^2 = \mathbf{j}^2 = \mathbf{k}^2 =  -1$,
\item $\mathbf{i} \mathbf{j} = - \mathbf{j} \mathbf{i} = \mathbf{k}$,
\item $\mathbf{j} \mathbf{k} = - \mathbf{k} \mathbf{j} = \mathbf{i}$,
\item $\mathbf{k} \mathbf{i} = - \mathbf{i} \mathbf{k} = \mathbf{j}$
\end{itemize}
the quaternion $q$ may be written as
$$
q = q_0 + q_1\mathbf{i} + q_2\mathbf{j} + q_3\mathbf{k}
$$
Hence, addition of quaternions is common
$$
p+q = (q_0 + q_0) + (p_1+q_1)\mathbf{i} + (p_2+q_2)\mathbf{j} + (p_3+q_3)\mathbf{k}
$$
while their product is \emph{not} commutative
\begin{equation}\label{eq:qProduct}
  \footnotesize
\left.\begin{array}{rcl}
p~q &=& {\color{white}+} \left(q_0 q_0 - (p_1 q_1 + p_2 q_2 + p_3 q_3) \right){\color{white}\mathbf{j}} + \left(p_0 q_1 + q_0 p_1 + (p_2 q_3 - q_2 p_3) \right)\mathbf{i} \\ 
&~&+ \left(p_0 q_2 + q_0 p_2 + (p_3 q_1 - q_3 p_1) \right)\mathbf{j} + \left(p_0 q_3 + q_0 p_3 + (p_1 q_2 - q_1 p_2) \right)\mathbf{k}
\end{array}\right.
\end{equation}
The norm of a quaternion $q$ is defined by means of its conjugate $q^*$
\begin{align*}
q^* &= q_0 - (q_1 \mathbf{i} + q_2 \mathbf{j} + q_3 \mathbf{k})\\
|q|^2 &= q~q^* = q^* ~ q
\end{align*}
Then, it follows that the inverse of $q$ is
$$
q^{-1} = \dfrac{q^*}{|q|^2}
$$

The imaginary part $q_1\mathbf{i} + q_2\mathbf{j} + q_3\mathbf{k}$ of a quaternion $q$ may be identified as a vector $\mathbf{q} \in \Re^3$
$$
q = q_0 + \mathbf{q}
$$
Using dot ($\cdot$) and cross ($\times$) products of vectors $\mathbf{p},\mathbf{q} \in \Re^3$, we may write \eqref{eq:qProduct} as
$$
p~q = \left(p_0 q_0 - \mathbf{p} \cdot \mathbf{q}\right) + \underline{\left(p_0 \mathbf{q} + q_0 \mathbf{p} + \mathbf{p} \times \mathbf{q} \right)}
$$
where the underlined part is the 3D vector representing the imaginary part of the product.
Scalar and vector products as well as vector addition are covariant with rotations of $\text{SO}(3)$, which implies that the imaginary part of a quaternion is also covariant.
It means that every rotation of the imaginary part $\mathbf{q}$ corresponds to an \emph{automorphism} of quaternions, i.e. a bijective mapping from quaternions to quaternions which preserves the structure of quaternions (i.e. their product).
We are going to identify this automorphism, which is related to quaternion product.

Assuming that the four components defining a quaternion correspond to coordinates of an euclidean space of fourth dimension, we consider unit quaternions $\hat{q}$ such that
$$
\hat{q}_0^2 + \hat{q}_1^2 + \hat{q}_2^2  + \hat{q}_3^2 = 1 = \hat{q}_0^2 + |\hat{\mathbf{q}}|^2
$$
A unit quaternion $\hat{q}$ may define right ($\mathcal{R}$) and left ($\mathcal{L}$) screws on a quaternion $p$ (whose norm may be different from 1), 
\begin{equation}\label{eq:screws}
\mathcal{R}_{\hat{q}}(p): p \mapsto \hat{q}~p, \hspace*{.5cm} \mathcal{L}_{\hat{q}}(p): p \mapsto p~\hat{q}
\end{equation}
which are automorphisms of quaternions\footnote{Actually, the definition of $\mathcal{L}$ will be modified to be an automorphism.}.

The screws do not alter the norm of $p$
$$
\begin{array}{rcl}
|p~\hat{q}|^2 &=& \left(p_0^2 \hat{q}_0^2  + (\mathbf{p} \cdot \hat{\mathbf{q}})^2  - 2 p_0 \hat{q}_0 \mathbf{p} \cdot \hat{\mathbf{q}} \right) +\\
& &	 \left(p_0^2 |\hat{\mathbf{q}}|^2 + \hat{q}_0^2 |\mathbf{p}|^2  + |\mathbf{p}|^2 |\hat{\mathbf{q}}|^2- (\mathbf{p} \cdot \hat{\mathbf{q}})^2 + 2 p_0 \hat{q}_0 \hat{\mathbf{q}} \cdot \mathbf{p} \right) \\
&=& (p_0^2 + |\mathbf{p}|^2) (\hat{q}_0^2 + |\hat{\mathbf{q}}|^2)\\
&=& |\hat{q}~p|^2\\
&=& |p|^2
\end{array}
$$
where we used the relationship $|\mathbf{a} \times \mathbf{b}|^2 = |\mathbf{a}|^2 |\mathbf{b}|^2 - (\mathbf{a} \cdot \mathbf{b})^2 $, with  $\mathbf{a}, \mathbf{b} \in \Re^3$.

Besides, if we consider the angle $\gamma$ between the fourth dimensional vectors corresponding to quaternions $p$ and $q$, it is the same between their image produced by any screws $\hat{r}$
$$
\begin{array}{rcl}
(\hat{r}~p) ~\odot_{\Re^4}~ (\hat{r}~q) &=&  \left(\hat{r}_0 p_0 - \hat{\mathbf{r}} \cdot \mathbf{p} + \hat{r}_0 \mathbf{p} + p_0 \hat{\mathbf{r}} + \hat{\mathbf{r}} \times \mathbf{p} \right) \\ && \odot_{\Re^4}  \left(\hat{r}_0 q_0 - \hat{\mathbf{r}} \cdot \mathbf{q} + \hat{r}_0 \mathbf{q} + q_0 \hat{\mathbf{r}} + \hat{\mathbf{r}} \times \mathbf{q}\right) \\
~&=& (\hat{r}_0^2+|\hat{\mathbf{r}}|^2) (p_0 q_0 + \mathbf{p} \cdot \mathbf{q}) \\
~&=& (p~\hat{r}) ~\odot_{\Re^4}~ (q~\hat{r}) \\
~&=& p ~\odot_{\Re^4} ~ q\\
~&=& \cos(\gamma) |p|~|q|
\end{array}
$$
where we used the following relationships
\begin{itemize}
\item $\mathbf{a} \cdot (\mathbf{b} \times \mathbf{c}) =  \mathbf{b} \cdot (\mathbf{c} \times \mathbf{a}) = \mathbf{c} \cdot (\mathbf{a} \times \mathbf{b})$
\item $(\mathbf{a} \times \mathbf{b}) \cdot (\mathbf{c} \times \mathbf{d}) = (\mathbf{a} \cdot \mathbf{c}) (\mathbf{b} \cdot \mathbf{d}) - (\mathbf{a} \cdot \mathbf{d}) (\mathbf{b} \cdot \mathbf{c})$
\end{itemize}
We understand that screws \eqref{eq:screws} correspond to some kind of rotations.
We are going to describe their properties starting with a well chosen unit quaternion $\hat{q} = \cos(\alpha) + \sin(\alpha) \mathbf{k}$.

The corresponding screws give
\begin{align*}
(p_0 \cos(\alpha) - p_3 \sin(\alpha)) &+ (p_1 \cos(\beta) - p_2 \sin(\beta))\mathbf{i} \\
+ (p_2 \cos(\beta) - p_1 \sin(\beta) )\mathbf{j} &+ (p_3 \cos(\alpha) + p_0 \sin(\alpha) )\mathbf{k}
\end{align*}
with $\beta=\pm \alpha $ for respectively $\mathcal{L}_{\hat{q}}(p)$ and $\mathcal{R}_{\hat{q}}(p)$.
For general values of $\beta$, the transformation may be seen as acting on a fourth dimensional vector, i.e. a matrix-vector product whose matrix is
$$
\begin{bmatrix}
\cos(\alpha) & 0 & 0 & -\sin(\alpha) \\
0 & \cos(\beta) & -\sin(\beta) & 0  \\
0 & \sin(\beta) & \cos(\beta) & 0 \\
\sin(\alpha) & 0 & 0 & \cos(\alpha) \\
\end{bmatrix}
$$
It is a \emph{compound} rotation of angles $\alpha$ in the $<1;\mathbf{k}>$-plane and $\beta$ in the $<\mathbf{i};\mathbf{j}>$-plane.
Those two planes are \emph{absolutely} orthogonal, which means that they have \emph{no} nonzero vector in common.
Then, the compound rotations act only on vectors in their corresponding plane.
Those planes actually define the invariant planes of the corresponding rotations.
As we said that every rotation of the imaginary part $\mathbf{q}$ of a quaternion $q$ corresponds to an automorphism of quaternions, our arbitrary choice $\hat{q} = \cos(\alpha) + \sin(\alpha) \mathbf{k}$ may be generalized to $\hat{r} = \cos(\alpha) + \sin(\alpha) \mathbf{v}$ with $\mathbf{v}\in\Re^3~ \text{s.t.}|\mathbf{v}|^2=1$.
The two invariant planes of $\mathcal{R}_{\hat{r}}(p), \mathcal{L}_{\hat{r}}(p)$ correspond to the one joining the real axis to $\mathbf{v}$ and the one that is normal to $\mathbf{v}$ contained within the region $<\mathbf{i};\mathbf{j};\mathbf{k}> = \Re^3\subset\Re^4$ (called \emph{imaginary prime}), since we have rotated the imaginary part of $q$ only. It means that the invariant planes are still absolutely orthogonal.

$\mathcal{R}_{\hat{r}}(p)~\forall \hat{r}$ is isomorphic to the product of unit quaternions denoted by $\hat{Q}$.
It means that the combination of $\mathcal{R}_{\hat{a}},\mathcal{R}_{\hat{b}}$ corresponds to $\mathcal{R}_{(\hat{b}\hat{a})}$, i.e. the right screw parameterized by the product of units quaternions $\hat{r}_b~\hat{r}_a$.
However, $\mathcal{L}_{\hat{r}}(p)$ is not isomorphic to the product of unit quaternions; it is not \emph{yet} an automorphism.
Indeed, the combination of $\mathcal{L}_{\hat{a}},\mathcal{L}_{\hat{b}}$ (in this order) has to act from right in the quaternion product.
To ensure a consistent combination, we have to modify the left screw definition, which is \emph{now} an automorphism
$$
\mathcal{L'}_{\hat{r}}(p): p \mapsto p \hat{r}^{-1}
$$
The combination 
$$\mathcal{L'}_{\hat{b} \hat{a}}(p): p \mapsto p \hat{a}^{-1}~\hat{b}^{-1}$$ is then consistent, since the left screws act on the left, in the correct order.

Owing to the associative multiplication of quaternions, combination of $\mathcal{R}_{\hat{r}}$ and $\mathcal{L'}_{\hat{r}}$ is simply
\begin{equation}\label{eq:automorphism}
\mathcal{A}_{\hat{r}}(p): p \mapsto \hat{r}~p~\hat{r}^{-1}
\end{equation}
Automorphisms of form \eqref{eq:automorphism} are \emph{inner} automorphisms.

Taking again the particular case $q = \cos(\alpha) + \sin(\alpha) \mathbf{k}$, $\mathcal{A}_{\hat{q}}$ may be expressed with products of the two matrices
{\scriptsize
\begin{align*}
&\underbrace{\begin{bmatrix}
\cos(\alpha) & 0 & 0 & \sin(\alpha) \\
0 & \cos(\alpha) & -\sin(\alpha) & 0  \\
0 & \sin(\alpha) & \cos(\alpha) & 0 \\
-\sin(\alpha) & 0 & 0 & \cos(\alpha) \\
\end{bmatrix}}_{\mathcal{L'}_{\hat{q}}}
\underbrace{\begin{bmatrix}
\cos(\alpha) & 0 & 0 & -\sin(\alpha) \\
0 & \cos(\alpha) & -\sin(\alpha) & 0  \\
0 & \sin(\alpha) & \cos(\alpha) & 0 \\
\sin(\alpha) & 0 & 0 & \cos(\alpha) \\
\end{bmatrix}}_{\mathcal{R}_{\hat{q}}} \\
=&\underbrace{\begin{bmatrix}
1 & 0 & 0 & 0 \\
0 & \cos(2\alpha) & -\sin(2\alpha) & 0\\
0 & \sin(2\alpha)  & \cos(2\alpha) &0 \\
0 & 0 & 0 & 1\\
\end{bmatrix}}_{\mathcal{A}_{\hat{q}}}
\end{align*}}
First, observe that the matrices may commute since the compound rotations act in two same invariant planes.
We notice that the two rotations in the $<1;\mathbf{k}>$-plane avoid each other, while in the $<\mathbf{i};\mathbf{j}>$-plane they add to each other.
Reminding that $\hat{q}$ is arbitrary and may be generalized to $\hat{r} = \cos(\alpha) + \mathbf{v} \sin(\alpha)$, we understand that $\mathcal{R}_{\hat{r}}$ is a compound rotation of angle $\alpha$ in the planes joining the real axis to $\mathbf{v}$ and the one which is normal to $\mathbf{v}$ contained in the imaginary prime, while $\mathcal{L'}_{\hat{r}}$ is respectively a rotation of $-\alpha$ and $\alpha$ in those latter planes. Hence, $\mathcal{A}_{\hat{r}}$ is a simple rotation of angle $2\alpha$ around $\mathbf{v}$ in the imaginary prime, which coincides with $\Re^3$.
We have identified the automorphism of quaternions which corresponds to a rotation of the imaginary part of a quaternion, leaving unchanged its real part: it is the inner automorphism \eqref{eq:automorphism}.
Those inner automorphisms provides us a \emph{homomorphic} mapping between $\hat{Q}$ and $\text{SO}(3)$, since two (opposite) unit quaternions correspond to a single rotation of $\text{SO}(3)$: $\pm \hat{r}$.
\begin{equation}\label{eq:q2so3}
\hat{Q} \xmapsto{2:1} \text{SO}(3)
\end{equation}
Obviously, the opposite is also true: any inner automorphisms of nonzero quaternions $q \in \Re^4$ represented by
\begin{equation}\label{eq:inner}
\mathcal{A}_{q}(p): p \mapsto q~p~q^{-1}
\end{equation}
is a rotation of the imaginary prime, since any nonzero $q\in \Re^4$ may be written as $q:=|q|\hat{q}$.

\subsection{Special unitary group $\text{SU}(2)$}

Let a quaternion $q=q_0 + q_1 \mathbf{i} + q_2 \mathbf{j} + q_3 \mathbf{k}$.
We define $u,v \in \mathbb{C}$
\begin{equation}\label{eq:q2c}
u := q_0 + q_3~i, ~ v := -(q_2 + q_1~i)
\end{equation}
such that they parameterize the complex matrix of the following form
\begin{equation}\label{eq:cm}
\begin{bmatrix}
u & -v^* \\
v & u^*
\end{bmatrix}
\end{equation}

The above relationship \eqref{eq:q2c} defines an isomorphism between the quaternions and the two-by-two skew-Hermitian matrices\footnote{$A\in \mathbb{C}^{2\times 2}$ s.t. $A_{ij}=-A_{ji}^*$} of form \eqref{eq:cm}, which are called \emph{quaternionic} matrices.
Let us consider two quaternions $p,q \in \Re^4$ and their corresponding quaternionic matrices $P,Q \in \mathbb{C}^{2 \times 2}$.
It is straightforward that the addition $P+Q$ corresponds $p+q$, and conversely.
The products $PQ$ and $p~q$ correspond to each other
\begin{align*}
\begin{bmatrix}
u_p & -v_p^* \\
v_p & u_p^*  \\
\end{bmatrix}
\begin{bmatrix}
u_q & -v_q^* \\
v_q & u_q^*  \\
\end{bmatrix} 
=
\begin{bmatrix}
u_p u_q - v_p^* v_q & - (u_p v_q^* + u_p^* u_q^*) \\
v_p u_q + u_p^* v_q & - v_p v_q^* + u_p^* u_q^*
\end{bmatrix}
\end{align*}
where
\begin{itemize}
\item {\footnotesize $u_p u_q - v_p^* v_q = {\color{white}-} \left(q_0 q_0 - (p_1 q_1 + p_2 q_2 + p_3 q_3)\right) + \left(p_0 q_3 + q_0 p_3 + (p_1 q_2 - q_1 p_2) \right) i$ }
\item {\footnotesize $v_p u_q + u_p^* v_q  = - \left(p_0 q_2 + q_0 p_2 + (p_3 q_1 - q_3 p_1) \right) -\left(p_0 q_1 + q_0 p_1 + (p_2 q_3 - q_2 p_3) \right)i$}
\end{itemize}
whose real and imaginary parts correspond to the quadruplet defining the product $p~q$.

If we apply \eqref{eq:q2c} to a unit quaternion $\hat{q} \in \hat{Q}$, the quaternionic \eqref{eq:cm} matrix has a determinant $\hat{u} \hat{u}^* + \hat{v} \hat{v}^*=1$.
This matrix form corresponds to the $\text{SU}(2)$ group.
Owing to \eqref{eq:q2c} for unit quaternions, we have
\begin{equation}\label{eq:su2rot}
\text{SU}(2) \overset{1:1}{\longleftrightarrow} \hat{Q} \xmapsto{2:1} \text{SO}(3)
\end{equation}

At this point, we could wonder why using \eqref{eq:q2c} to define \eqref{eq:cm}.
If we use the common definition $u,v \in \mathbb{C}$
$$
u = q_0 + q_1~i, ~ v = q_2 + q_3~i
$$
the product $p~q$ would correspond to the product of quaternionic matrices $QP$ (instead of $PQ$).

Let us come back to screws $\mathcal{L}$ and $\mathcal{R}$.
We consider $\hat{q},\hat{p} \in \hat{Q}$ and $r \in \Re^4$, with their respective complex coordinates $(\hat{u}_q;\hat{v}_q),(\hat{u}_p;\hat{v}_p) \in \text{SU}(2)$ and $(u_r;v_r) \in \mathbb{C}^2$ defined by \eqref{eq:q2c}.
We know those screws correspond to rotations of the fourth dimensional vector corresponding to $\Re^4$, which may be represented by a matrix-vector product.
The corresponding quaternionic representation of $\mathcal{R}_{\hat{r}}(p)$ is
\begin{equation}\label{eq:affineTransformation}
\begin{bmatrix}
u_p \\ v_p
\end{bmatrix}
\mapsto
\begin{bmatrix}
\hat{u}_{\hat{r}} & -\hat{v}_{\hat{r}}^* \\ \hat{v}_{\hat{r}} & \hat{u}_{\hat{r}}^*
\end{bmatrix}
\begin{bmatrix}
u_p \\ v_p
\end{bmatrix}
=
\begin{bmatrix}
\hat{u}_{\hat{r}} u_p - \hat{v}_{\hat{r}}^* v_p \\ \hat{v}_{\hat{r}} u_p + \hat{u}_{\hat{r}}^* v_p
\end{bmatrix}
\end{equation}
whose real and imaginary parts correspond to the four components of $\hat{p}~r \in \Re^4$.
The transformation 	\eqref{eq:affineTransformation} is a \emph{complex affine transformation}.
However, there is no such affine transformation (acting on $(u_p;v_p)$) corresponding to the left screw $\mathcal{L'}_{\hat{r}}$.
The corresponding quaternionic representation is
$$
\begin{bmatrix}
u_p & -v_p^*
\end{bmatrix}
\mapsto
\begin{bmatrix}
u_p & -v_p^*
\end{bmatrix}
\begin{bmatrix}
\hat{u}_{\hat{r}} & -\hat{v}_{\hat{r}}^* \\ \hat{v}_{\hat{r}} & \hat{u}_{\hat{r}}^*
\end{bmatrix}^{-1}
=
\begin{bmatrix}
\hat{u}_{\hat{r}}^* u_p + \hat{v}_{\hat{r}} v_p^* & \hat{v}_{\hat{r}}^* u_p - \hat{u}_{\hat{r}} v_p^*
\end{bmatrix}
$$
which is a complex affine transformation on $(u_p;-v_p^*)$.

%%%%%%%%%%%%%%%%%%
\section{Frame Field Space}

In the last section, we have seen that a frame can be represented by a rotation of the cartesian frame.
However, there are 24 distinct rotations (i.e. 48 unit quaternions) giving a unique frame.
We need functions which give a unique set of values for each frame, in other words a same set of values for the 48 unit quaternions generating a given frame.

We are going to analyze how some rotational groups  act on $(u;v) \in \mathbb{C}^2$ as affine transformations.
This analysis will give three invariant forms per group, giving a set of values corresponding to coordinates of a surface embedded in $\mathbb{C}^3$.
We eventually need to redefine $(u;v)$ from a quaternion, in order to properly parameterize the frames from those invariant forms.

\subsection{Finite groups of quaternions}

We are mainly interested in the octahedral group, which shares the symmetries of a frame.
Nevertheless, we need to define it from two smaller finite groups of quaternions, which are the vierer and binary tetrahedral groups.

\FloatBarrier

\emph{The vierer group} $\mathcal{V} \subset \hat{Q}$, only consists of 4 rotations that are of angle $\pi$ around the axes of the cartesian frame including the identity
$$
\mathcal{V} = \{\pm 1; \pm \mathbf{i}; \pm \mathbf{j}; \pm \mathbf{k}\}
$$

\FloatBarrier
\emph{The binary tetrahedral} group $\mathcal{T} \subset \hat{Q}$ is composed of the 12 rotations that leave unchanged the orientation of a regular tetrahedron whose 4 vertices are located at $ (1;1;1), (1;-1;-1), (-1;1;-1), (-1;-1;1)$ and of its dual whose vertices have respectively opposite components of the primal.
$$
\mathcal{T} = \mathcal{V} \oplus \left\{ \dfrac{1}{2}\left( \pm 1 \pm \mathbf{i} \pm \mathbf{j} \pm \mathbf{k}\right)\right\}
$$

\begin{figure}[!ht]
  \begin{center}
    \subfloat[Primal tetrahedron.]{\includegraphics[width=.45\linewidth]{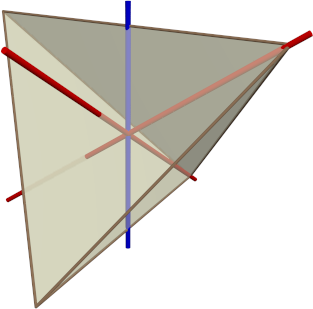}}\hspace*{.5cm}
    \subfloat[Dual tetrahedron.]{\includegraphics[width=.45\linewidth]{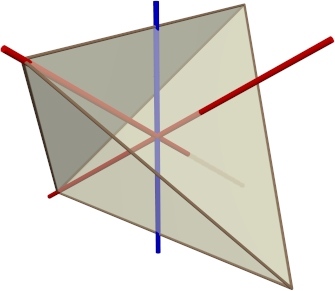}}
    \caption{Three axes of rotations generating the binary tetrahedral group. Blue and red axes correspond respectively to rotations of $\pi$ and $\frac{2\pi}{3}$.}
    \label{fig:binaryTetrahedral}
  \end{center}
\end{figure}
\FloatBarrier

\emph{The binary octahedral} group $\mathcal{O} \subset \hat{Q}$ has 24 rotations that preserve the orientation of an octahedron whose vertices are $(\pm 1;0;0),(0;\pm 1;0),\\ (0;0;\pm 1)$ and of the dual cube whose the centers of each face correspond to the vertices of the primal octahedron.
\begin{align*}
\mathcal{O} = \mathcal{T} \oplus   \left\{ \dfrac{1}{\sqrt{2}}\left(\pm 1 \pm \mathbf{i}\right); \dfrac{1}{\sqrt{2}}\left(\pm 1 \pm \mathbf{j}\right); \dfrac{1}{\sqrt{2}}\left(\pm 1 \pm \mathbf{k}\right);\right.  \\ \left. \dfrac{1}{\sqrt{2}}\left(\pm \mathbf{i} \pm \mathbf{j}\right); \dfrac{1}{\sqrt{2}}\left(\pm \mathbf{j} \pm \mathbf{k}\right); \dfrac{1}{\sqrt{2}}\left(\pm \mathbf{i} \pm \mathbf{k}\right) \right\}
\end{align*}

\begin{figure}[!ht]
  \begin{center}
    \subfloat[Hexahedron.]{\includegraphics[width=.45\linewidth]{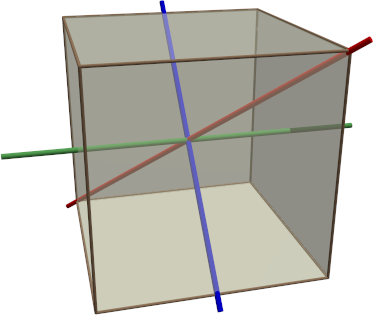}}\hspace*{.5cm}
    \subfloat[Octahedron.]{\includegraphics[width=.45\linewidth]{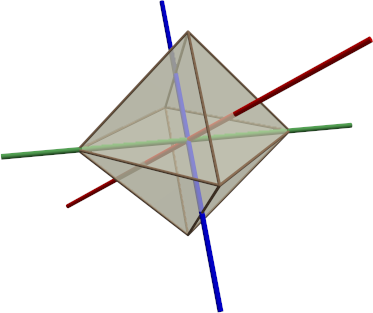}}
    \caption{Three axes of rotations generating the binary octahedral group. Blue, red and green axes correspond respectively to rotations of $\pi, \frac{2\pi}{3}$ and $\frac{\pi}{2}$.}
    \label{fig:binaryTetrahedral}
  \end{center}
\end{figure}

\FloatBarrier
Each finite group may be generated by three unit quaternions
{\footnotesize\begin{equation}\label{eq:Qpqr}
\hat{a} = \cos\left(\frac{\pi}{p}\right) + \mathbf{a} \sin\left(\frac{\pi}{p}\right),~ \hat{b} = \cos\left(\frac{\pi}{q}\right) + \mathbf{b} \sin\left(\frac{\pi}{q}\right),~ \hat{c} = \cos\left(\frac{\pi}{2}\right) + \mathbf{c} \sin\left(\frac{\pi}{2}\right)
\end{equation}}
with $\mathbf{a},\mathbf{b},\mathbf{c}$ the corresponding axes of rotation.
Those unit quaternions are such that they satisfy the relationships
\begin{equation}\label{eq:pqr}
\hat{a}^p = \hat{b}^q = \hat{c}^2 = \hat{a} \hat{b} \hat{c} = -1
\end{equation}
which is depicted by the triplet $(p;q;2)$; it corresponds to the powers of \eqref{eq:pqr}.

The triplet $(2;2;2)$ corresponds to the vierer group $\mathcal{V}$, with (e.g.) $\mathbf{a}=\mathbf{i}, \mathbf{b}=\mathbf{j}, \mathbf{c}=\mathbf{k}$.
The triplet $(3;3;2)$ corresponds to the binary tetrahedral group $\mathcal{T}$, with (e.g.) $\mathbf{a}=\frac{1}{\sqrt{3}} (\mathbf{i}+\mathbf{j}+\mathbf{k}), \mathbf{b}=\frac{1}{\sqrt{3}} (-\mathbf{i}+\mathbf{j}+\mathbf{k}), \mathbf{c}=\mathbf{k}$.
The triplet $(4;3;2)$ corresponds to the binary octahedral group $\mathcal{O}$, with (e.g.) $\mathbf{a}= \mathbf{i}, \mathbf{b}=\frac{1}{\sqrt{3}} (\mathbf{i}+\mathbf{j}+\mathbf{k}), \mathbf{c}=\frac{1}{\sqrt{2}} (\mathbf{j}-\mathbf{k})$.

\subsection{Model surface}

We look for a set of functions $f_k: \mathbb{C}^2 \mapsto \mathbb{C}$ such that 
\begin{equation}\label{eq:invariantAim}
f_k\left(A_i \begin{bmatrix}
u \\ v
\end{bmatrix}\right)
= 
f\left( \begin{bmatrix}
u \\ v
\end{bmatrix}\right),~\forall A_i \in \mathcal{V}, \mathcal{T}, \text{ xor } \mathcal{O}
\end{equation}
where $A_i$ is a quaternionic matrix corresponding to an element of the finite group $(p;q;2)$, i.e.
$$
A_i^p=-I, A_i^q=-I, \text{ xor } A_i^2 = -I
$$
with $I$ the two-by-two identity matrix.

We know that a finite group $(p;q;2)$ may be generated by three unit quaternions written as \eqref{eq:Qpqr}.
Hence, if \eqref{eq:invariantAim} is satisfied for such three unit quaternions, it is satisfied by all quaternions of the finite group. 
It means we have three unit vectors, i.e. $6$ real degrees of freedom.
It implies that the set of functions is composed of three functions $(f_0;f_1;f_2) \in \mathbb{C}^3$.
However, the last equality of \eqref{eq:pqr} is not satisfy for all quaternions written as  \eqref{eq:Qpqr}.
It is a constraint, that will appears as a polynomial relationship between $f_0,f_1$ and $f_2$.
\begin{equation}\label{eq:modelSurface}
f_2^2 = \mathcal{P}\left(f_0;f_1\right)
\end{equation}
where $\mathcal{P}$ is a bivariate polynomial.
This relationship defines a surface embedded in $\mathbb{C}^3$.
It is parameterized by $(u;v)\in \mathbb{C}^2$.
Besides, we observe that there exists an affine transformation on $(u;v)$ which does not alter $f_0(u;v)$ and $f_1(u;v)$, but which gives an opposite sign to $f_2(u;v)$.
This latter transformation does not belong to the current finite group $(p;q;2)$.

Such a surface \eqref{eq:modelSurface} is a \emph{model surface}, which defines the rotations of a quotient group:  $\hat{Q}/\mathcal{V}, \hat{Q}/\mathcal{T}~\text{or}~\hat{Q}/\mathcal{O}$ (abusing of the isomorphism between $\text{SU}(2)$ and $\hat{Q}$).
The functions $f_k$ parameterizing it are homogeneous polynomial in $(u;v)\in\mathbb{C}^2$.

\subsection{Invariant forms}

We are going to derive the invariant forms relative to each finite group.

\subsubsection*{The vierer group}
The vierer group $\mathcal{V}$ has the identity rotation, which may be represented by the following affine transformation
$$
\begin{bmatrix}
u \\ v
\end{bmatrix}
\mapsto
\begin{bmatrix}
-1 & 0 \\ 0 & -1
\end{bmatrix}
\begin{bmatrix}
u \\ v
\end{bmatrix}
=
\begin{bmatrix}
-u \\ -v
\end{bmatrix}
$$
We find that $uv$, $u^2$ and $v^2$ are invariant under such a transformation.
However, they do change for half turns around axes $\mathbf{i},\mathbf{j},\mathbf{k}$: in the case of $\mathbf{i}$
$$
\begin{bmatrix}
u \\ v
\end{bmatrix}
\mapsto
\begin{bmatrix}
0 & i \\ i & 0
\end{bmatrix}
\begin{bmatrix}
u \\ v
\end{bmatrix}
=
\begin{bmatrix}
i~v \\ i~u
\end{bmatrix}
$$
which respectively gives $-uv$, $-v^2$ and $-u^2$.
We observe that we have to square every expression and then to combine the two last ones and thus to make a third expression from the three initial ones: $(uv)^2$, $u^4+v^4$ and $uv(u^2+v^2)(u^2-v^2)$.
Those latter expressions are also invariant for half turns around $\mathbf{j}$ and $\mathbf{k}$.
They are thus invariant forms of $\mathcal{V}$
\begin{equation}\label{eq:viererForms}
\begin{array}{l}
\left\{\begin{array}{rcl}

f_0(u;v) &=& (uv)^2\\ f_1(u;v) &=& u^4 + v^4\\ f_2(u;v) &=& uv(u^4-v^4)\end{array}\right.
\\
\begin{array}{rcl}
\text{s.t.} ~ f_2^2(u;v) &=& f_0(u;v) \left(f_1^2(u;v) - 4 f_0^2(u;v)\right)
\end{array}
\end{array}
\end{equation}
As claimed, there is an affine transformation that does not belong to $\mathcal{V}$ that does not alter $f_0(u;v),f_1(u;v)$ but which changes the sign of $f_2(u;v)$
$$
\begin{bmatrix}
u \\ v
\end{bmatrix}
\mapsto
\begin{bmatrix}
1 & 0 \\  0 & -1
\end{bmatrix}
\begin{bmatrix}
u \\ v
\end{bmatrix}
=
\begin{bmatrix}
u \\ -v
\end{bmatrix}
$$

\subsubsection*{The binary tetrahedral group}
The binary tetrahedral group $\mathcal{T}$ contains $\mathcal{V}$.
Since a quaternion group is generated by three unit quaternions,  if we build invariant forms for $\mathcal{T}$ from those of $\mathcal{V}$, we just need to have them invariant for two unit quaternions of $\mathcal{T}$ (since it will be invariant to unit quaternions of $\mathcal{V}$) such that $(3;3;2)$.
We consider first the rotation $\frac{1}{2}(1+\mathbf{i}+\mathbf{j}+\mathbf{k})$ described by the following transformation
$$
\begin{bmatrix}
u \\ v
\end{bmatrix}
\mapsto
\dfrac{1}{2}\begin{bmatrix}
1+i & i-1 \\  1+i & 1-i
\end{bmatrix}
\begin{bmatrix}
u \\ v
\end{bmatrix}
=
\dfrac{1}{2}(1+i)\begin{bmatrix}
u+i~v \\ u-i~v
\end{bmatrix}
$$
We first notice that  $\tilde{f}_0 = f_1 + 2\sqrt{3}i ~f_0 = u^4 + v^4 + 2\sqrt{3}i (uv)^2,~\tilde{f}_1 = f_1 - 2\sqrt{3}i ~f_0  = u^4 + v^4 - 2\sqrt{3}i (uv)^2,~\tilde{f}_2=f_2$ are also invariant forms of $\mathcal{V}$ s.t. $\tilde{f}_2^2 = \frac{\sqrt{3}i}{36} (\tilde{f}_1^3 - \tilde{f}_0^3)$. We may see that
$$\begin{array}{rcl}
\tilde{f}_0(u+i~v; u-i~v) &=& (u+iv)^4 + (u-iv)^4 + 2\sqrt{3}i (u^2+v^2)^2 \\
&=& 2(1+\sqrt{3}i) (u^4+v^4) + 4(\sqrt{3}i-3) (uv)^2\\
&=& 2(1+\sqrt{3}i) (u^4+v^4 + 2\sqrt{3}i (uv)^2)\\
&=&  4 \exp\left(\frac{\pi i}{3}\right)\tilde{f}_0(u;v)
\end{array}
$$
and that
$$\begin{array}{rcl}
\tilde{f}_1(u+i~v; u-i~v) &=& (u+iv)^4 + (u-iv)^4 - 2\sqrt{3}i (u^2+v^2)^2 \\
&=& 2(1-\sqrt{3}i) (u^4+v^4) - 4(3+\sqrt{3}i) (uv)^2\\
&=& 2(1-\sqrt{3}i) (u^4+v^4 - 2\sqrt{3}i (uv)^2)\\
&=&  4 \exp(-\frac{\pi i}{3}) \tilde{f}_0(u;v)
\end{array}
$$
Knowing that $\left(\dfrac{1+i}{2}\right)^4=-\dfrac{1}{4}$, we notice that
\begin{itemize}
\item {\scriptsize $ \tilde{f}_0\left(\frac{1}{2}(1+i)(u+iv);\frac{1}{2}(1+i)(u-iv)\right)\tilde{f}_1\left(\frac{1}{2}(1+i)(u+iv);\frac{1}{2}(1+i)(u-iv)\right) = (\tilde{f}_0\tilde{f}_1)(u;v)$}
\item {\scriptsize $\tilde{f}_0^3\left(\frac{1}{2}(1+i)(u+iv);\frac{1}{2}(1+i)(u-iv)\right) = \tilde{f}_0^3(u;v)$}
\item {\scriptsize $\tilde{f}_1^3\left(\frac{1}{2}(1+i)(u+iv);\frac{1}{2}(1+i)(u-iv)\right) = \tilde{f}_1^3(u;v)$}
\item {\scriptsize $\tilde{f}_2\left(\frac{1}{2}(1+i)(u+iv);\frac{1}{2}(1+i)(u-iv)\right) = \tilde{f}_2(u;v)$}
\end{itemize}
Since $\tilde{f}_i$ are homogeneous polynomials in $(u;v)$, we get the same results for $\frac{1}{2}(1-i+j-k)$.
It is then possible to build invariant forms from the above relationships.
Indeed, the product $\tilde{f}_0 \tilde{f}_1$ and $\tilde{f}_2$ are invariant for elements of $\mathcal{V}$, $\frac{1}{2}(1+i+j+k)$ and $\frac{1}{2}(1-i+j-k)$ which satisfy $(3;3;2)$.
Those two expressions are thus invariant to all elements of $\mathcal{T}$.
Remembering that $\tilde{f}_2^2 =  \frac{\sqrt{3}i}{36} (\tilde{f}_1^3 - \tilde{f}_0^3) = (uv(u^4-v^4))^2$ is a linear combination of $\tilde{f}_0^3$ and $\tilde{f}_1^3$, which are invariant to $\mathcal{T}$, we end up with the following invariant forms for $\mathcal{T}$
\begin{equation}\label{eq:tetrahedralForms}
\begin{array}{l}
\left\{\begin{array}{rcl}

g_0(u;v) &=& uv(u^4-v^4) = \tilde{f}_2(u;v) \\ g_1(u;v) &=& u^8 + v^8 + 14 (uv)^4 = \tilde{f}_0(u;v) \tilde{f}_1(u;v) \\ g_2(u;v) &=& u^{12} + v^{12} - 33 (uv)^4 (u^4+v^4) = \tilde{f}_0^3(u;v) + \tilde{f}_1^3(u;v)\end{array}\right.
\\
\begin{array}{rcl}
\text{s.t.} ~ g_2^2(u;v) &=& g_1^3(u;v) - 108 g_0^4(u;v)
\end{array}
\end{array}
\end{equation}
The affine transformation
\begin{equation}\label{eq:affineSign}
\begin{bmatrix}
u \\ v
\end{bmatrix}
\mapsto
\begin{bmatrix}
\exp\left(\frac{\pi i}{4}\right) & 0 \\  0 & -\exp\left(-\frac{\pi i}{4}\right)
\end{bmatrix}
\begin{bmatrix}
u \\ v
\end{bmatrix}
=
\begin{bmatrix}
\exp\left(\frac{\pi i}{4}\right) u \\ - \exp\left(-\frac{\pi i}{4}\right) v
\end{bmatrix}
\end{equation}
leaves $g_0$ and $g_1$ unchanged while $g_2$ has its sign changed.

\subsubsection*{The binary octahedral group}
The binary octahedral group $\mathcal{O}$ containing $\mathcal{T}$, which contains $\mathcal{V}$.
Since $\mathcal{O}$ may be generated by three unit quaternions such that $(4;3;2)$, we just have to build invariant forms from $\mathcal{T}$ that are invariant to a rotation consisting in a quarter turn around one of the axes $\mathbf{i}, \mathbf{j}, \mathbf{k}$ (or any combinations of two of them).
In the case of a quarter turn around $\mathbf{i}$, it corresponds to the following affine transformation
$$
\begin{bmatrix}
u \\ v
\end{bmatrix}
\mapsto
\begin{bmatrix}
\exp\left(\frac{\pi i}{4}\right) & 0 \\  0 & \exp\left(-\frac{\pi i}{4}\right)
\end{bmatrix}
\begin{bmatrix}
u \\ v
\end{bmatrix}
=
\begin{bmatrix}
\exp\left(\frac{\pi i}{4}\right) u \\ \exp\left(\frac{-\pi i}{4}\right) v
\end{bmatrix}
$$
which obviously leaves unchanged $g_1(u;v)$, but it changes the sign of $g_0(u;v)$ and $g_2(u;v)$.
It means that $g_0^2$ and the product $g_0 g_1$ are invariant under the above affine transformation.
Therefore, we get the following invariant forms for $\mathcal{O}$ 
\begin{equation}\label{eq:octahedralForms}
\begin{array}{l}
\left\{\begin{array}{rcl}

h_0(u;v) &=& u^8 + v^8 + 14(uv)^4 = g_1(u;v) \\ h_1(u;v) &=& (uv(u^4-v^4)^2 = g_0^2(u;v) \\ h_2(u;v) &=& uv(u^4-v^4)(u^{12}-33u^4v^4(u^4+v^4)+v^{12})=  (g_0 g_2)(u;v)\end{array}\right.
\\
\begin{array}{rcl}
\text{s.t.} ~ h_2^2(u;v) &=& h_1(u;v) \left(h_0^3(u;v) - 108 h_1^2(u;v)\right)
\end{array}
\end{array}
\end{equation}
$h_0(u;v)$ and $h_1(u;v)$ are invariant under the affine transformation \eqref{eq:affineSign}, while $h_2(u;v)$ takes the opposite sign.

\subsection{Going back to frames}
We have identified invariant forms for \emph{right} screws taking their value in either $\mathcal{V}, \mathcal{T}$ or $\mathcal{O}$.
In the case of $\mathcal{O}$, it is expressed as
$$
\hat{o}_i~\hat{r} = \hat{q}_i
$$
with $\hat{o}_i \in \mathcal{O}, 0\le i < 48$ and $\hat{r} \in \hat{Q}$.
Hence, the 48 rotations of the cartesian frame $\mathbf{\hat{f}} = \{\pm \mathbf{i}; \pm \mathbf{j} ; \pm \mathbf{k}\}$ giving another 3D frame are described by  $\hat{q}_i$
$$
(\hat{o}_i~\hat{r})~\mathbf{\hat{f}}~(\hat{o}_i~\hat{r})^{-1}
$$
However, these rotations do not describe a 3D frame.
The rotations $\hat{o}_i, \hat{r}$ are performed in the wrong order.
The 48 rotations $\hat{o}_i$ have to be applied first on $\mathbf{\hat{f}}$, then only the rotation $\hat{r}$ giving the new 3D frame has to be applied.

To reverse the order, we have to conjugate the quaternions $\hat{q}_i$ corresponding to $(h_0;h_1;h_2)$.
Indeed,
\begin{equation}\label{eq:reversing}
(\hat{o}_i~\hat{r})^*~\mathbf{\hat{f}}~((\hat{o}_i~\hat{r})^*)^{-1} = \hat{r}^*~\hat{o}_i^* ~ \mathbf{\hat{f}} ~ \hat{o}_i~\hat{r}
\end{equation}
Since $\hat{o}_i^* \in \mathcal{O} ~ \forall i$, \eqref{eq:reversing} applies in the correct order the rotations, which produces a new 3D frame corresponding to a rotation $\hat{r}$ of $\mathbf{\hat{f}}$ up to a symmetry.
Therefore, we have to update \eqref{eq:q2c} in order to take account of the conjugation of $\hat{q}$
\begin{equation}\label{eq:conjugating2reverse}
\left.\begin{array}{rcl}
u &=& q_0 - q_1~i \\
v &=& q_2 + q_3~i
\end{array}\right\}
\end{equation}

\section{Numerical Insights}

A set of 48 $(\hat{u};\hat{v}) \in \text{SU}(2)$ which preserves the orientation of a given cube is a \emph{groupset}.
We have shown that all elements of a groupset are mapped onto the same complex valued coordinates $(h_0;h_1;h_2)\in \mathbb{C}^3$ of the model surface \eqref{eq:octahedralForms}.
It is possible to do the reverse way: from any coordinate of the model surface, the corresponding groupset can be identified.
The groupset may be composed of nonunit quaternions; however, all the 48 quaternions have the same norm owing to the fact that the corresponding affine transformations do not alter the norm.

We are going to rewind the derivations, starting with $h_0,h_1,h_2$ and going through $g_0,g_1,g_2$ and $\tilde{f}_0,\tilde{f}_1$ in order to get $u,v$.
Owing to \eqref{eq:octahedralForms}, we get
$$
\begin{array}{rcl}
g_1 &=& h_0\\
g_0 &=& (-1)^{k_0} \sqrt{h_1}\\
g_2 &=& 
\left\{\begin{array}{lr}
(-1)^{k_0}\frac{h_2}{\sqrt{h_1}} & ,~\text{if}~ |h_1| \neq 0\\
(-1)^{k_0} \sqrt{h_0^3} &,~\text{otherwise}
\end{array}\right.
\end{array}
$$
with $k_0=\{0;1\}$. There are thus two possibilities.
Then, we obtain with \eqref{eq:tetrahedralForms}
$$
\begin{array}{rcl}
\tilde{f}_0 &=& \exp\left(i\frac{2\pi}{3}k_1\right)  \sqrt[\leftroot{-1}\uproot{2} 3]{\dfrac{2 g_2 + 12\sqrt{3}i g_0^2}{2}}\\
\tilde{f}_1 &=& 
\left\{\begin{array}{lr}
\frac{g_1}{\tilde{f}_0} &,~\text{if}~|\tilde{f}_0| \neq 0\\
\exp\left(i\frac{2\pi}{3}k_1\right)  \sqrt[\leftroot{-1}\uproot{2} 3]{\dfrac{2 g_2 - 12\sqrt{3}i g_0^2}{2}} &,~\text{otherwise}
\end{array}\right.
\end{array}
$$
with $k_1=\{0,1,2\}$. There are thus three possibilities.
Knowing that $\tilde{f}_{0,1}(u;v) = u^4 \pm 4\sqrt{3}i (uv)^2 + v^4$ (respectively), we may identify 
$$
(uv)^2 = \dfrac{\tilde{f}_0-\tilde{f}_1}{4\sqrt{3}i}
$$
\subsubsection*{$\bullet$ $\tilde{f}_0\neq\tilde{f_1}$}
We can compute $v^4$ from the following quadratic polynomial
$$
48 (v^4)^2 - 24 (\tilde{f}_0 + \tilde{f}_1) v^4 - (\tilde{f}_0-\tilde{f}_1)^2 = 0
$$
whose roots are
$$
v = i^{k_3} \sqrt[\leftroot{-1}\uproot{2} 4]{ \frac{\tilde{f}_{0}}{4} + \frac{\tilde{f}_{1}}{4} + (-1)^{k_2} \frac{\sqrt{3} \sqrt{\tilde{f}_{0}^{2} + \tilde{f}_{0} \tilde{f}_{1} + \tilde{f}_{1}^{2}}}{6} }
$$
with $k_2= \{0,1\}, k_3 = \{0,1,2,3\}$.
In that case, there are eight possibilities for $v$, which each gives
$$
u = (-1)^{k_4} \sqrt{v}
$$
with $k_4=\{0,1\}$.
We end up with $96$ possible $(u;v)$, which is twice what we expected.
It is due to the choice of the sign of the square root in the expression $g_0 = (-1)^{k_0} \sqrt{h_1}$.
Indeed, in one case we choose the wrong sign for $g_0$, which gives the wrong sign to $g_2$.
A wrong choice gives the groupset defined by $(h_0;h_1;{\color{red}-}h_2)$, instead of $(h_0;h_1;h_2)$.
A wrong choice may be easily mapped onto an element of $(h_0;h_1;h_2)$, by using the affine transformation \eqref{eq:affineSign}.
That latter transformation is performed in practice, since we only need one quaternion to rotate the cartesian frame onto the underlined 3D frame.

\subsubsection*{$\bullet$ $\tilde{f}_0=\tilde{f_1}$}
Otherwise if $uv$ is zero, it means that $u=0$ or $v=0$, giving
$$
\left\{\begin{array}{rcl}
  u &=& k_2 ~ i^{k_3}  \sqrt[\leftroot{-1}\uproot{2} 4]{\frac{\tilde{f}_0+\tilde{f}_1}{2}}\\
  v &=& (1-k_2)~ i^{k_3} \sqrt[\leftroot{-1}\uproot{2} 4]{\frac{\tilde{f}_0+\tilde{f}_1}{2}}
        \end{array}\right.
$$
        with $k_2= \{0,1\}, k_3 = \{0,1,2,3\}$.
There are then eight possibilities of $u$ and $v$.
We end up with 48 possible $(u;v)$.
We do not encounter the latter issue, because if $u$ or $v$ is zero, $h_2=0$.

%%%%%%%%%%%%%%%%%%%%%
\subsection{Euclidean distance $\mathbb{C}^3$}\label{sec:euclideanDistance}

We analyze here the distance between 3D frames $\mathbf{\tilde{f}}$ and the cartesian frame $\mathbf{\hat{f}}$.
The aim is to check if $(h_0;h_1;h_2)\in\mathbb{C}^3$ such that \eqref{eq:octahedralForms} may be used to define an euclidean distance between 3D frames.

\begin{figure}
\begin{center}
\subfloat[3D frames are frame rotated around one axis.]{\includegraphics[scale=.49]{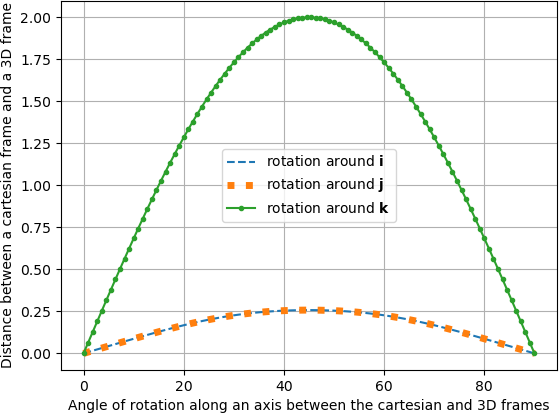}\label{sub:distanceAxis}}\\
\subfloat[Distances between $10^4$ random 3D frames and the cartesian frame.]{\includegraphics[scale=.49]{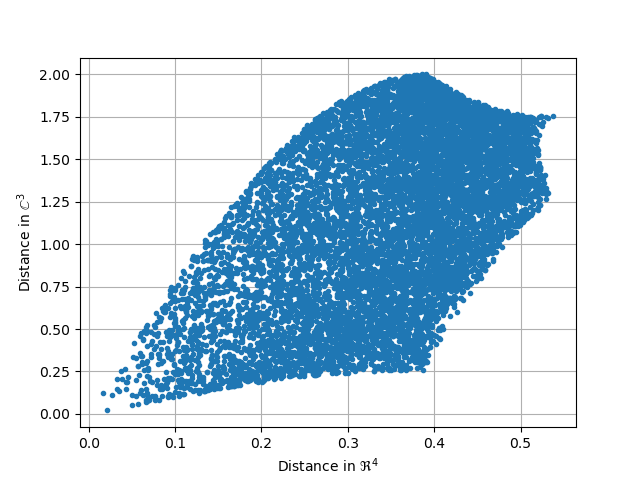}\label{sub:C3vsR4}}
\caption{Distance between 3D frames and the cartesian frame. The distances are the euclidean ones of $\mathbb{C}^3$ (to compare $(h_0;h_1;h_2)$) and of $\Re^4$ (to compare $(q_0;q_1;q_2;q_3)$.)}
\end{center}
\end{figure}

First, we produce frames rotated around a single axis,
$$
\mathbf{\tilde{f}} = \left[ \cos\left(\dfrac{\alpha}{2}\right) + \mathbf{v} \sin\left(\dfrac{\alpha}{2}\right)  \right] ~ \mathbf{\hat{f}} ~ \left[ \cos\left(\dfrac{\alpha}{2}\right) - \mathbf{v} \sin\left(\dfrac{\alpha}{2}\right)  \right]
$$
with $\mathbf{v} \in \{\mathbf{i};\mathbf{j};\mathbf{k}\}$ and $\alpha \in \left[0;\frac{\pi}{2}\right]$.
We then compute their distance to the cartesian frame
\begin{equation}\label{eq:frameDistance}
d_{\mathbb{C}^3}(\mathbf{\tilde{f}};\mathbf{\hat{f}}) = \sqrt{ dh_0 dh_0^* + dh_1 dh_1^* + dh_2 dh_2^* }
\end{equation}
with $dh_i = \tilde{h}_i - \hat{h}_i$ the difference of the i-th components of the triplet \eqref{eq:octahedralForms}.
Fig. \ref{sub:distanceAxis} shows that frames rotated around $\mathbf{k}$ appear to be further than the ones rotated around $\mathbf{i},\mathbf{j}$.
Obviously, the former ones should be as far than the latter ones.
\eqref{eq:frameDistance} is not isotropic.

Let us compare \eqref{eq:frameDistance} with the euclidean distance in $\Re^4$, i.e. the shortest distance between any of 48 unit quaternions giving $\mathbf{\tilde{f}}$ from $\mathbf{\hat{f}}$ and the unit quaternion $1$
$$
d_{\Re^4}(\mathbf{\tilde{f}};\mathbf{\hat{f}}) = \underset{i}{\min} \sqrt{(\hat{q}_{0(i)}-1)^2 + \hat{q}_{1(i)}^2 + \hat{q}_{2(i)}^2 + \hat{q}_{3(i)}^2 + \hat{q}_{4(i)}^2}
$$
This actually defines a consistent distance between frames.
The frames $\mathbf{\tilde{f}}$ are produced from random unit quaternions $\hat{q}$.
Fig. \ref{sub:C3vsR4} shows that the distances do not correspond at all.
Frames having triplet \eqref{eq:octahedralForms} close to the one of $\mathbf{\hat{f}}$ may be either close or far of the cartesian frame.

We may conclude \eqref{eq:frameDistance} is not suitable to define a distance between frames.
It means that averaging sets of values $(h_0;h_1;h_2)$ in $\mathbb{C}^3$ is inconvenient.
Indeed, two sets that are close according to \eqref{eq:frameDistance} could give a frame which is not the average of the two frames corresponding to the two sets.

\subsection{Around an axis}

We analyze the behavior of $h_0,h_1,h_2$ of 3D frames having one of their axes in common.
We define this axis with $\mathbf{n} =(n_x;n_y;n_z) \in \Re^3$ such that $n_x^2 + n_y^2 + n_z^2 = 1$.
Tanks to the automorphism of quaternions, we consider
$$
\hat{q}_z = \cos\left(\dfrac{\alpha}{2}\right) + \sin\left(\dfrac{\alpha}{2}\right) \mathbf{k}
$$
parameterizing rotations around $\mathbf{k}$ of angle $\alpha$.
Afterwards, we rotate $\mathbf{k}$ onto $\mathbf{n}=(n_x;n_y;n_z) \in \Re^3$ by means of $\hat{q}_{z \mapsto n} $.
We thus aim
$$
\hat{q}_n = \hat{q}_{z \mapsto n} ~ \hat{q}_z
$$

The unit quaternion $\hat{q}_{z \mapsto n}$ may be described by a rotation of $\gamma$ around a unit vector $\mathbf{v}$.
It gives
\begin{align*}
\cos(\gamma) &= \mathbf{k} \cdot \mathbf{n} = n_z\\
\sin(\gamma) &= |\mathbf{k} \times \mathbf{n}| = \sqrt{n_x^2 + n_y^2} = \sqrt{1-n_z}\sqrt{1+n_z}\\
\mathbf{v} &= \dfrac{\mathbf{k} \times \mathbf{n}}{|\mathbf{k} \times \mathbf{n}|} = \dfrac{(-n_y;n_x;0)}{\sqrt{1-n_z}\sqrt{1+n_z}}
\end{align*}
Using trigonometric identities for $\cos$ and $\sin$ of $\frac{\gamma}{2}$, we get
$$
\hat{q}_{z \mapsto n} = \dfrac{1}{\sqrt{2}} \left( \sqrt{1+n_z} + \dfrac{1}{\sqrt{1+n_z}} ( -n_y \mathbf{i} + n_x \mathbf{j} )  \right)
$$
We then get
$$
\begin{array}{rcl}
  \hat{q}_n &=& \dfrac{1}{\sqrt{2}} \left( \sqrt{1+n_z} \cos\left(\frac{\alpha}{2}\right) +  \dfrac{1}{\sqrt{1+n_z}} [ n_x \sin\left(\frac{\alpha}{2}\right) - n_y \cos\left(\frac{\alpha}{2}\right) ]\mathbf{i} \right. \\
  && \left.+ \dfrac{1}{\sqrt{1+n_z}} [ n_x \cos\left(\frac{\alpha}{2}\right) + n_y \sin\left(\frac{\alpha}{2}\right)]\mathbf{j} + \sqrt{1+n_z} \sin\left(\frac{\alpha}{2}\right) \mathbf{k} \right)
\end{array}
$$
Owing to \eqref{eq:conjugating2reverse}, we have the following $\text{SU}(2)$ representation
\begin{align*}
\hat{u}_n &= \sqrt{\dfrac{1+n_z}{2}} \exp\left(-i\frac{\alpha}{2}\right) \\
\hat{v}_n &= \dfrac{n_x - i n_y}{\sqrt{2(1+n_z)}} \exp\left(i\frac{\alpha}{2}\right)
\end{align*}
It then gives the invariant forms
\begin{equation}\label{eq:axisInvariantForm}
  \scriptsize
\begin{array}{rcl}
h_{0[n]}(\alpha) &=& \left(\dfrac{n_x - i n_y}{2}\right)^4 ~ \left( \left[\dfrac{n_x-in_y}{n_z+1} \exp(i\alpha)\right]^4 + 14 +  \left[\dfrac{n_z+1}{n_x-in_y} \exp(-i\alpha)\right]^4  \right) \\
h_{1[n]}(\alpha) &=& \left(\dfrac{n_x - i n_y}{2}\right)^6 ~ \left( \left[\dfrac{n_x - i n_y}{n_z + 1} \exp(i\alpha)\right]^2 - \left[\dfrac{n_z + 1}{n_x - i n_y} \exp(-i\alpha)\right]^2 \right)^2 \\
h_{2[n]}(\alpha) &=& - \left(\dfrac{n_x - i n_y}{2}\right)^9 ~ \left( \left[\dfrac{n_x - i n_y}{n_z + 1} \exp(i\alpha)\right]^2 - \left[\dfrac{n_z + 1}{n_x - i n_y} \exp(-i\alpha)\right]^2  \right) \\ && \left(  \left[\dfrac{n_x - i n_y}{n_z + 1} \exp(i\alpha)\right]^{6} + \left[\dfrac{n_z + 1}{n_x - i n_y} \exp(-i\alpha)\right]^{6} \right. \\ && \left. - 33 \left\{  \left[\dfrac{n_x - i n_y}{n_z + 1} \exp(i\alpha)\right]^{2} + \left[\dfrac{n_z + 1}{n_x - i n_y} \exp(-i\alpha)\right]^{2} \right\}  \right)
\end{array}
\end{equation}
Let
$$
\left.\begin{array}{rcl}
w(\alpha) &:=& \dfrac{n_x-in_y}{n_z+1} \exp(i\alpha)\\
a &:=& \dfrac{n_x - i n_y}{2}
\end{array}\right\}
$$
which simplifies the writing of \eqref{eq:axisInvariantForm}
\begin{equation}\label{eq:axisInvariantFormSimplified}
\begin{array}{|rcl|}
\hline
h_0(w) &=& a^4 (w^4 + w^{-4} + 14)\\
\hline
h_1(w) &=& a^6 (w^{2}-w^{-2})^2 \\ &=& a^6 (w^4 + w^{-4} - 2)\\
\hline
h_2(w) &=& -a^9 (w^2-w^{-2}) (w^6 + w^{-6} - 33 (w^2+ w^{-2})) \\
&=& -a^9 (w^8 - w^{-8} - 34 (w^4 - w^{-4})) \\
&=& -a^9 (w^4 - w^{-4}) (w^4 + w^{-4} - 34)\\
\hline
\end{array}
\end{equation}
We  notice that
\begin{align*}
a^2 h_0(w) - h_1(w) &= 16 a^6 \\
-a^3 (w^4 - w^{-4}) (h_1(w)-32) &= h_2(w)
\end{align*}
giving a linear relationship between $h_0$ and $h_1$ parameterized by two first components of the axis $\mathbf{n}$, and updating the model surface equation.

\begin{figure}
\begin{center}
\includegraphics[scale=.5]{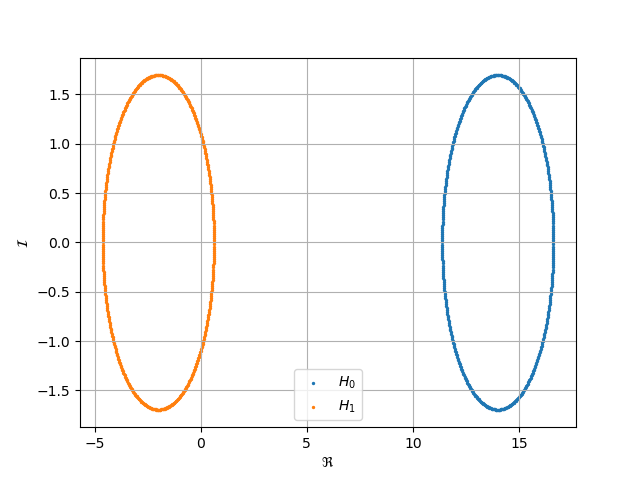}
\caption{Scattering of $H_0,H_1$ in the complex plane $\mathbb{C}$ for 500 three-dimensional frames rotated around an axis $\mathbf{n}$ which corresponds to one of the axes of those frames.}
\label{fig:axisEllipsis}
\end{center}
\end{figure}

Let us define
$$
H_0 := \dfrac{h_0}{a^4},~H_1 := \dfrac{h_1}{a^6}
$$
and write 
$$
\dfrac{n_x - i n_y}{n_z+1} = \dfrac{1-n_z}{1+n_z} \exp(i~\theta)
$$
in polar coordinates with $\theta=\arctan\left(\dfrac{-n_y}{n_x}\right)$.
We then let $t = 4(\alpha+\theta)$ to obtain
$$
\begin{array}{rcl}
H_0(t) &=& p \exp(i t) + q \exp(-i t) + 14 \\
H_1(t) &=& p \exp(i t) + q \exp(-i t) - 2 \\
\end{array}
$$
with $p=\dfrac{1}{q}=\left(\dfrac{1-n_z}{1+n_z}\right)^4$.
We notice $H_0, H_1$ have expressions corresponding to ellipses in the complex plane, Fig. \ref{fig:axisEllipsis}.
It follows that $h_0,h_1$ are ellipses which are scaled and rotated.
The corresponding relationships could be used to set boundary conditions on frames in order to align one of their axes with the normal of the boundary of a volume.

\subsection{Ensuring $\text{SU}(2)$ from model surface}\label{sec:ensuring}

We are interested in unit quaternions.
We know that a coordinate of the model surface provides a groupset whose quaternions have the same norm.
We could compute a quaternion corresponding to a coordinate by using the above procedure.
However, it is possible to derive a simpler way to check if a given coordinate $(h_0;h_1;h_2)$ laying on the model surface corresponds to $(\hat{u};\hat{v})\in \text{SU}(2)$.

We do not need to compute explicitly a quaternion.
We just need to know the value of some powers of $u$ and $v$.
We remind that
$$
|q|^2_{\Re^4} = |u|^2_{\mathbb{C}} + |v|^2_{\mathbb{C}}
$$
We are going to derive a way to compute those values $u^n$ and $v^m$.

First, we notice that applying affine transformations (as right screws) from the octahedral group to $f_0(u;v)=(uv)^2$\footnote{We remind that $f_0$ is an invariant form of the vierer group.} produces\footnote{$f_0\left(A_i \cdot \begin{bmatrix}u \\ v \end{bmatrix}\right)=\dots$} the following subset
$$
(uv)^2,~ -\dfrac{1}{4}(u^2+v^2),~ \dfrac{1}{4}(u^2-v^2)^2
$$

The three elements of this subset may be written as the roots of polynomial of third order $(x-(uv)^2)(x+(u^2+v^2)/4)(x-(u^2-v^2)/4)=0$.
Scaling and expanding the latter expression allows to identify $h_0$ and $h_1$ as coefficients of the polynomial expression
\begin{equation}\label{eq:polynomialUV2}
16 x^3 - h_0 x + h_1 = 0
\end{equation}

Let $x^\star :=(uv)^2$ for an arbitrary root value of the above polynomial; choosing an other root value accounts of choosing an other quaternion of the same groupset.

We then compute either $u^8$ and $v^8$ as the two roots of the quadratic polynomial $(y-u^8)(y- v^8)=0$.
Expanding that latter polynomial in terms of $h_0$ and $x^\star$, we get
\begin{equation}\label{eq:polynomialU8}
y^2 - (h_0-14{x^\star}^2)y + {x^\star}^4 = 0
\end{equation}

Then the norm of the corresponding $(u;v)$ is given by
$$
|(u;v)|^2 = |u^8|^{\frac{1}{4}} + |v^8|^{\frac{1}{4}}
$$
Knowing that the discriminant of \eqref{eq:polynomialU8} is
$$
\Delta^2 = h_0^2-28h_0 {x^\star}^2 + 192 {x^\star}^4
$$
we have the following norm
$$
|(u;v)|^2 = \left| \dfrac{h_0 - 14 {x^\star}^2 + \Delta}{2} \right|^{\frac{1}{4}} + \left| \dfrac{h_0 - 14 {x^\star}^2 - \Delta}{2} \right|^{\frac{1}{4}} 
$$
which is independent of the choice of root for $x^\star$ and $\Delta$, since modifying a choice corresponds to choosing an other quaternion of the same groupset.

\subsection{Attempted numerical schemes}

As mentionned in \S \ref{sec:euclideanDistance}, the use of the eulidean distance related to $\mathbb{C}^3$ for measuring the distance between frames from their triplet $(h_0;h_1;h_2)$  is unconsistent.
Yet, we have tried different numerical schemes based on a finite element approach.
Even if those have been unsuccessful due to this latter statement, we describe them.

The schemes are based on a tetrahedral mesh that discretizes the region $R$ of interest.
A Lagrange $\mathcal{P}^1$ approximation is built from that mesh.

%%%%%%%%%%%%%%%%%%%%%%%%%%%%%%%%%%
\subsubsection{Linear formulation}

We consider that the frames laying on the boundary $\partial R$ are given.
In pratice we compute the corresponding crossfield\citep{Beaufort2017Jan}, and then identify a rotation sending the cartesian frame to the frame.
From that rotation, we have a corresponding complex pair $(u;v)\in\mathbb{C}^2$.
We eventually have $(h_0(u;v);h_1(u;v);h_2(u;v))$ all over the boundary $\partial R$.

We assume that $h_i$ is function of $(x;y;z) \in R \subset \Re^3, ~\forall i$.
We aim to get smooth values of $h_i$ within $R$, which corresponds to minimize their Dirichlet energy
$$
\underset{(h_0;h_1;h_2)}{\min}\int_R  \left| \nabla h_0(x;y;z) \right|_{\mathbb{C}}^2 + \left| \nabla h_1(x;y;z) \right|_{\mathbb{C}}^2 + \left| \nabla h_2(x;y;z) \right|_{\mathbb{C}}^2 ~  dx dy dz
$$
with $|\bullet|_{\mathbb{C}}$ the usual complex norm.
The weak finite formulation is then given by
$$
\sum_{i=0}^2 \sum_j \int_R  \nabla \phi_j \cdot \nabla \phi_k ~ dx dy dz ~ h_{i,j} = 0, \forall k
$$
with $h_{i,j}$ the nodal value of $h_i$ in node $(x_j;y_j;z_j)$.

We get a linear system with three complex unknows by vertex (node).
The solution is not projected onto the model surface.

%%%%%%%%%%%%%%%%%%%%%%%%%%%%%%%%%%
\subsubsection{Collocation method}\label{sec:collocation}

We parameterize the rotations using the Euler angles $(\alpha,\beta,\gamma)$ respectively around $\mathbf{k},\mathbf{j},\mathbf{k}$, i.e. the following matrix belonging to $\text{SO}(3)$
{\tiny \begin{equation}\label{eq:eulerKJK}
  \begin{bmatrix}
    - \sin{\left(\alpha \right)} \sin{\left(\gamma \right)} + \cos{\left(\alpha \right)} \cos{\left(\beta \right)} \cos{\left(\gamma \right)} & -\sin{\left(\alpha \right)} \cos{\left(\gamma \right)} - \sin{\left(\gamma \right)} \cos{\left(\alpha \right)} \cos{\left(\beta \right)} & \sin{\left(\beta \right)} \cos{\left(\alpha \right)}\\
     \sin{\left(\alpha \right)} \cos{\left(\beta \right)} \cos{\left(\gamma \right)} + \sin{\left(\gamma \right)} \cos{\left(\alpha \right)} & - \sin{\left(\alpha \right)} \sin{\left(\gamma \right)} \cos{\left(\beta \right)} + \cos{\left(\alpha \right)} \cos{\left(\gamma \right)} & \sin{\left(\alpha \right)} \sin{\left(\beta \right)}\\
    -\sin{\left(\beta \right)} \cos{\left(\gamma \right)} & \sin{\left(\beta \right)} \sin{\left(\gamma \right)} & \cos{\left(\beta \right)}
  \end{bmatrix}
\end{equation}}
whose the columns correspond to the 3 directions of a frame.
We know that \eqref{eq:eulerKJK} is equivalent to two opposite quaternions,
{\tiny $$
\pm \left( \cos\left(\dfrac{\beta}{2}\right)\cos\left(\dfrac{\alpha+\gamma}{2}\right) ; \sin\left(\dfrac{\beta}{2}\right)\sin\left(\dfrac{\gamma-\alpha}{2}\right) ; \sin\left(\dfrac{\beta}{2}\right)\cos\left(\dfrac{\gamma-\alpha}{2}\right) ; \cos\left(\dfrac{\beta}{2}\right)\sin\left(\dfrac{\alpha+\gamma}{2}\right) \right)
$$}

From this quaternion, we use the relationships \eqref{eq:conjugating2reverse} and \eqref{eq:octahedralForms} in order to get the octahedral forms parameterized with $(\alpha;\beta;\gamma)$.
Thanks to this parameterization, we may express the minimization of the Dirichlet energy based on the Euler angles
\begin{equation}\label{eq:collocation}
\underset{(\alpha_j;\beta_j;\gamma_j)}{\min} \sum_{i=0}^2 \sum_j \int_R |h_{i}(\alpha_j;\beta_j;\gamma_j) \nabla \phi_j(x;y;z)|_\mathbb{C}^2 ~ dxdydz
\end{equation}
which is a nonlinear optimization problem.

We solve \eqref{eq:collocation} using a Newton's method.
The required derivatives are computed by means of the chain rule based on the Wirtinger calculus\citep{wirtinger1927formalen}.
Owing to $|f|_\mathbb{C}^2=f \cdot f^*$, we rewrite \eqref{eq:collocation}
$$
\underset{(\alpha_l;\beta_l;\gamma_l)}{\min} \sum_{i=0}^2  \int_R  \left(\sum_j h_{i,j} \nabla \phi_j \right) \cdot \left(\sum_k h_{i,k}^* \nabla \phi_k\right) ~ dx dy dz
$$
with $h_{i,l}=h_i(\alpha_l;\beta_l;\gamma_l)$.
The corresponding Euler-Lagrange equations are thus
$$
\displaystyle{\sum_{i=0}^2 \sum_j \int_R \left( \dfrac{\partial h_{i,l}}{\partial \bullet_l} h_{i,j}^*  + h_{i,j}  \dfrac{\partial h_{i,l}^*}{\partial \bullet_l}  \right)  \nabla \phi_l\cdot \nabla \phi_j ~ dx dy dz} = 0, \forall l
$$
for $\bullet_l \in \{\alpha_l;\beta_l;\gamma_l\}$, with
$$
\dfrac{\partial h_{i,l}}{\partial \bullet_l} = \dfrac{\partial h_{i,l}}{\partial u} \dfrac{\partial u}{\partial \bullet_l} + \dfrac{\partial h_{i,l}}{\partial v} \dfrac{\partial v}{\partial \bullet_l}
$$

We eventually need the corresponding hessian\footnote{for the entry $(k,l)$}
{\fontsize{5}{6}\selectfont$$
\displaystyle{\sum_{i=0}^2 \int_R  \left( \dfrac{\partial h_{i,l}}{\partial \bullet_l}\dfrac{\partial h_{i,k}^*}{\partial \bullet_k} + \dfrac{\partial h_{i,k}}{\partial \bullet_k}\dfrac{\partial h_{i,l}^*}{\partial \bullet_l} + \sum_j \delta_{kl} \left( \dfrac{\partial^2 h_{i,l}}{\partial \bullet_l \partial \star_l} h_{i,j}^* + h_{i,j} \dfrac{\partial^2 h_{i,l}^*}{\partial \bullet_l \partial \star_l}  \right)\right) ~ \nabla \phi_l \cdot \nabla \phi_k ~ dx dy dz}
$$}
for $\bullet_l,\star_l \in \{\alpha_l;\beta_l;\gamma_l\}$ with
{\fontsize{5}{6}\selectfont$$
\dfrac{\partial h_{i,l}}{\partial \bullet_l \partial \star_l} = \dfrac{\partial^2 h_{i,l}}{\partial u^2} \dfrac{\partial u}{\partial \bullet_l}\dfrac{\partial u}{\partial \star_l} + \dfrac{\partial^2 h_{i,l}}{\partial v^2} \dfrac{\partial v}{\partial \bullet_l}\dfrac{\partial v}{\partial \star_l} + \dfrac{\partial h_{i,l}}{\partial u} \dfrac{\partial^2 u}{\partial \bullet \partial \star} + \dfrac{\partial h_{i,l}}{\partial v} \dfrac{\partial^2 v}{\partial \bullet \partial \star} + \dfrac{\partial^2 h_{i,l}}{\partial u \partial v} \left( \dfrac{\partial u}{\partial \bullet_l}\dfrac{\partial v}{\partial \star_l} + \dfrac{\partial v}{\partial \bullet_l}\dfrac{\partial u}{\partial \star_l} \right)
$$}

The boundary conditions are set by imposing that each frame on $\partial R$ shares a direction with the outward normal $\mathbf{n}=(n_x;n_y;n_z)$.
We consider that the last column of \eqref{eq:eulerKJK} corresponds to $\mathbf{n}$, which implies
$$
\begin{array}{rcl}
  \alpha &=& \left\{ \begin{array}{lc}
                       \arctan\left(\dfrac{n_y}{n_z}\right),&\text{if}~n_x\neq 0 \neq n_y \\
                       \dfrac{\pi}{2},&\text{if}~n_x=0,n_y \neq 0 \\
                       0,&\text{otherwise}
                     \end{array}\right.\\
 \beta &=& \pm \arccos(n_z), \text{the sign is determined with}~n_x~\text{or}~n_y
\end{array}
$$

%%%%%%%%%%%%%%%%%%%%%%%%%%%%%%%%%%%%%%
\subsubsection{Metric}

Fig. \ref{sub:C3vsR4} shows that the euclidean distance between the invariant forms is irrelevant, compared to the euclidean distance bewteen the quaternion (which is relevant).
We aim to measure the variation among the quaternions $\delta q$ from the variation of their invariant form $\delta h$.
A linear approximation is to built a metric $M$.
We start with the jacobian
$$
J = 
\begin{bmatrix}
  \dfrac{\partial h_0}{\partial u} & \dfrac{\partial h_0}{\partial v} \\ & \\
  \dfrac{\partial h_1}{\partial u} & \dfrac{\partial h_1}{\partial v} \\ & \\
  \dfrac{\partial h_2}{\partial u} & \dfrac{\partial h_2}{\partial v} 
\end{bmatrix}
$$
which measures the variation $\delta h$ from $\delta q$.
We thus need to compute a pseudo inverse $J^+$ of $J$.
Then, we can build the aimed metric $M$
$$
M = J^+ \cdot (J^+)^\dagger
$$
where $A^\dagger$ is the transposed conjugate (i.e. hermitian) of $A$.

Therefore, the Dirichlet energy corresponds to
\begin{equation}\label{eq:dirichletMetric}
\int_R \underline{\nabla h} \cdot M \cdot \underline{\nabla h}^\dagger  ~ dx dy dz
\end{equation}
where $\underline{\nabla h} = \begin{bmatrix} \nabla h_0 & \nabla h_1 & \nabla h_2 \end{bmatrix}$.

We can evaluate \eqref{eq:dirichletMetric} by means of a collocation method, as we have done in \S \ref{sec:collocation}.
The minimization is a bit trickier, since the metric $M$ is implicit, i.e. we cannot express it analytically from $(u;v)$, neither from $(\alpha;\beta;\gamma)$.
The derivatives have to be computed numerically.

However, \eqref{eq:dirichletMetric} remains an approximation which is only valid for points $(h_0;h_1;h_2)$ close to each other (laying in the neighborhood of the same tangent space).
The use of the metric is not sufficient to get a consistent distance of frames from their invariant forms.

\section{Conclusion}

We essentially gave a new parameterization of 3D frame fields, involving only 3 complex values related by an implicit equation describing a variety.
The SU(2) parameterization is based on \citep{du1964homographies}, with a slight modiftication about the isomorphism between the special unitary group and the unit quaternions.
We derived the invariant forms without involving homographies; we used the fact that a finite group of quaternions can be defined by three unit quaternions.
We showed how to get the quaternions from a coordinate of the variety.
Through numerical experiments, we noticed that the euclidean distance between 2 coordinates of the variety is not consistent with the distance between the corresponding 3D frames.
We derived the relationship of components for 3D frames sharing an even direction: the two first ones describe an ellipsis in their respective complex plane.
We showed how to ensure that a coordinate gives a unit quaternion.
Finally, three attempted numerical schemes have been described; they do not compute properly a frame field because of the inconsistent use of the euclidean norm.

\section*{Acknowledgements}

The present study was carried out in the framework of the project ``Large Scale
Simulation of Waves in Complex Media'', which is funded by the Communauté
Française de Belgique under contract ARC WAVES 15/19-03.
The authors thank Professor Jean-Pierre Tignol who has enlightened them about invariant theory.
\S\ref{sec:ensuring} is based on a Professor Tignol's note.

\bibliography{su2}

%-------------------------------------------------------------------------------------------------------------------------------------------------------------
\end{document}